\documentclass{amsart}[12pt]
\usepackage{amssymb,amsmath}
\usepackage{enumerate}

\usepackage[english]{babel}
\usepackage{color}
\usepackage{graphicx}

\newcommand {\gm} {\gamma}

\newcommand {\al} {\alpha}
\newcommand {\bt} {\beta}

\newcommand {\pr} {\prime}
\newcommand {\mc} {\mathcal}

\newcommand {\eqv} {\Longleftrightarrow}

\newtheorem{prob}{Problem}[section]

\theoremstyle{definition}
\newtheorem{rem}{Remark}[section]
\newtheorem{df}{Definition}[section]

\title{A contemporary inversive geometry \\solution of the three-circle problem of Steiner}
\begin{document}
\date{August 10, 2025}

\begin{abstract}
We present a concise self-contained inversive geometry solution of the three-circle problem of Steiner of constructing a circle that intersects each of the three given circles at one of the three given angles.
\end{abstract}

\author{AZIZKHON AZIZOV}
\address{The National University of Uzbekistan\\ Tashkent \\ Uzbekistan}
\email{a.azizov@nuu.uz}

\author{SEMYON LITVINOV}
\address{The Pennsylvania State University \\ 76 University Drive \\ Hazleton, PA 18202, USA}
\email{snl2@psu.edu}

\maketitle
\section {Introduction}

In the $3^{\text{rd}}$ century BC, Apollonius of Perga posted and solved the famous problem of constructing a circle tangent to three given circles, stated as Problem A in Subsection 5.1 of this article -- see also, for example, \cite{do, co,cg}. (Figuratively speaking, a geometric object is called {\it constructible} if it can be drawn using given geometric objects with the help of a straightedge and a compass only). However, Apollonius' solution the problem has not survived.

\vskip3pt
The significance of this problem can hardly be overstated -- no wonder the title of the Presidential address delivered in 1967 by the great geometer of the $20^{\text{th}}$ century H.S.M. Coxeter, at the joint meeting with the Mathematical Association of America, to the Canadian Mathematical Congress was ``The Problem of Apollonius'' At the beginning of his presentation Coxeter stated, ``As Sir Thomas Heath remarks \cite[p.\,182]{he}, this problem `has exercised the ingenuity of many distinguished geometers, including Vieta and Newton.'" (The results presented in this address were later published in \cite{co}.) Of note,
Vieta found a solution of the problem of Apollonius by considering limiting cases, where any of the three given circles can be reduced to a point or blown up into a line, which may have been the Apollonius' original approach.

\vskip3pt
This article is devoted to solving a generalized Apollonius problem of constructing a circle that intersects each of the three given circles at one of the three given  angles, from the interval $[0^0,90^0]$.

\vskip3pt
Only after we had presented our solution \cite{al}, it was brought to our attention that this ``three-circle problem" was already claimed long ago in \cite{st} by Jakob Steiner, the great geometer of the $19^{\text{th}}$ century -- see \cite[pp.\,468-469, Problem (b)]{am}, also \cite{lo}. On the bottom of page 162 in \cite{st}, Steiner states that he has solved this problem and that the solution will appear in a future monograph. However, it seems almost certain that that monograph was never published.

\vskip5pt
Notes section of this old Monthly \cite[p.\,469]{am} attributes the first published solution of this problem to F. Neumann in 1825 \cite{ne}. That Monthly note also contains quite a few related sources. For example, the book \cite{la} developes material at great length ``to {\it describe} a circle which shall cut three given circles at given angles" (pp.\,239-241), while \cite{coo} uses coaxiability properties of circles to solve ``the problem of Steiner," namely, ``to {\it construct} a circle meeting three given circles at three given angles."  Furthermore, we were also informed about the master's degree thesis \cite{pa}, where a detailed, 30-plus page, solution of the problem of Steiner that utilizes inversion is presented, which most likely was never published.

\vskip5pt
Our solution of the three-circle problem of Steiner -- stated as Problem S in Subsection 3.1 -- seems to be different from all that we have been able to find. We use inversive geometry -- briefly outlined in Section 4 of this article -- with the main tool being constructibility of a circle, the inversion with respect to which transforms two non-intersecting circles into concentric ones -- see Subsection 4.2. To solve the problem, we first assume, in Subsection 5.1, that a circle to be constructed exists and perform an inversion\, $\mc I$ of the four circles that transforms two of the three given circles into concentric ones. Then we find a formula for a possible radius of the image -- under the inversion $\mc I$ -- of the circle in question and formulas for possible distances from its center to the centers of the images of two of the given circles. Finally, using these three formulas, we {\it construct}, in Subsection 5.2, a circle $\mc C^\pr$ that does intersect each of the images of the three given ones at one of the three given angles. As the square of an inversion is the identity map, $\mc I(\mc C^\pr)$ is a circle that intersects each of the three given circles at one of the three given angles, finalizing the solution. Note that the case when each of the  three prescribed angles is of $90$ degrees is treated separately.

\vskip5pt
There is a way to approach the problem of Steiner by using coordinates to obtain a system of equations for the center and radius of the desired circle. But, just doing that and no more does not address the general existence or constructibility of this circle -- it simply allows for easy computation and graphing of that circle using software assuming its existence. Some analytic work has been done along these lines -- see for example \cite{se}. However, as the authors state in Preface to \cite{cg}, ``The mathematics curriculum in the secondary school normally includes a single one-year course in geometry," which is usually the student's ``sole exposure to the subject. In contrast, the mathematically minded student has the opportunity of studying elementary algebra, intermediate algebra, and even advanced algebra. It is natural, therefore, to expect a bias in favor of algebra and against geometry. Moreover, misguided enthusiasts lead the student to believe that geometry is ``outside the main stream of mathematics" and that analysis or set theory should supersede it."

\vskip5pt
Furtermore, the flavor of an analytic treatment would not be consonant with the original role of the problem as an application of the advanced concepts and techniques of synthetic Euclidean geometry that were developed in the late 18th and 19th centuries, such as inversive geometry. We have chosen to keep that spirit alive by presenting a solution to the three-circle problem with prescribed angles that uses these elegant synthetic geometric techniques to show not only that required circle exists, but also how it can be constructed.

\vskip5pt
Finally, aware of the expository challenge, we set up the explicit goal of producing an account that is self-contained by today's standards. The older solutions that we have seen are brief and elegant, but are given after lengthy developments of nontrivial synthetic geometric facts and techniques. In order to read them, one must do a lot of backtracking to gather up the necessary pieces before putting them together.

\section{Preliminaries}
\subsection{The (acute) angle between two circles}

Let $A$ be a point of intersection of two circles as in Figure 1. Then the acute angle at which these circles intersect is defined as the acute angle $\al$ between the tangent lines $t_1$ and $t_2$ to the circles at the point $A$. Clearly this angle coincides with the acute angle between the lines that pass through $O_1$ and $A$ and $O_2$ and $A$:
\vskip5pt
\begin{center}
\includegraphics[width=0.6\textwidth]{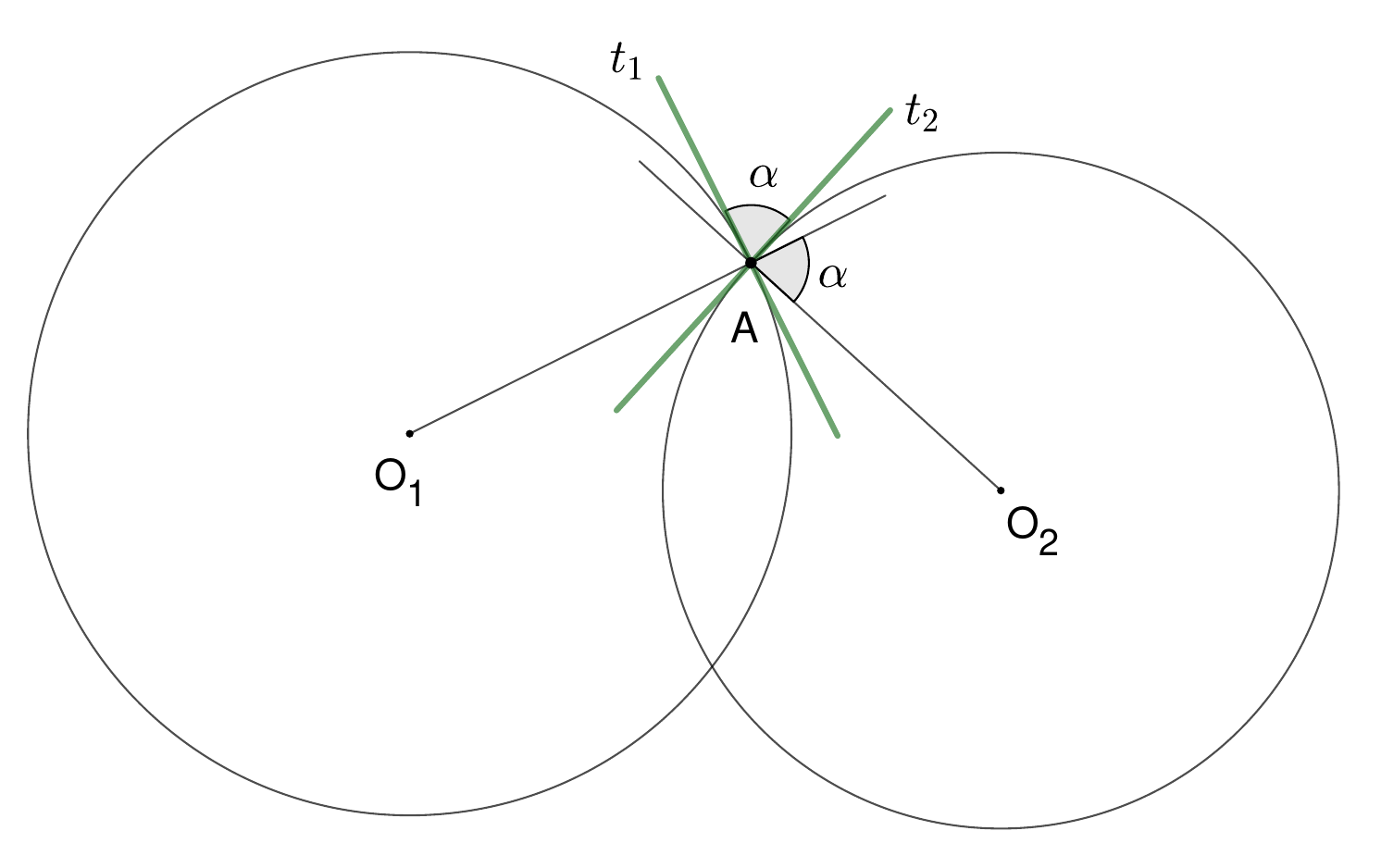}
\end{center}

\vspace{-5pt}
\centerline{Figure 1}

\subsection{Constructibility}

First off, we presume that the reader is familiar with the following compass-straightedge constructions:
\begin{enumerate}[(a)]
\item
Finding the mid-point of a given interval.
\item
Drawing the line through a given point orthogonal to a given line.
\item
Drawing the line through a given point parallel to a given line.
\end{enumerate}

\vskip5pt
Let us now try to define the notion of {\it constructibility} of an interval:

\begin{df}
Let $a_1, a_2,\dots, a_n$ be drawn intervals. Denote their lengths by the same letters. Let $f$ be a real-valued function of several real variables. An interval of length
\[
a=f(a_1, a_2,\dots, a_n)
\]
is called {\it constructible} if it can be drawn, using the intervals $a_1,\dots,a_n$, with the help of a straightedge and a compass only. In such a case, the number $a$ is also called constructible.
\end{df}

\noindent
A comprehensive account on constructibility can be found, for example, in \cite{cr}. We shall present just two basic constructions -- where we use those in (a)\,-\,(c) above without explicit references -- that we will need to attain the goal of this article:

\begin{prob}\label{prob1}
Construct the circle that passes through three non-collinear points $A$, $B$, and $C$.
\end{prob}
\noindent
{\bf Solution:} Utilizing constructions (a) and (b) above, draw perpendicular bisectors of the intervals $[AB]$ and $[BC]$.
The point of intersection of these bisectors $O$ is the center of the circle through $A$, $B$, and $C$, and we can draw this circle with a compass since its radius equals $|OA|$.

\begin{prob}\label{prob3}
Given intervals $a$, $b$, and $c$, construct an interval of length $x=\displaystyle\frac{ab}c$.
\end{prob}
\noindent
{\bf Solution:} Referring to Figure 2, draw segments $[OB]$ and $[OC]$ and the point $A\in[OC]$ such that $|OA|=a$, $|OC|=c$, and $|OB|=b$. Next, construct the line through $A$ parallel to $BC$, and let $X$ be the point of intersection of this line with $OB$. Then, denoting $|OX|=x$, we have
\[
\frac{|OX|}{|OA|}=\frac{|OB|}{|OC|}\ \ \eqv \ \ \frac xa=\frac bc\ \ \eqv \ \ x=\frac{ab}c.
\]
Therefore, the interval $[CX]$ has the desired length.

\begin{center}
\includegraphics[width=0.4\textwidth]{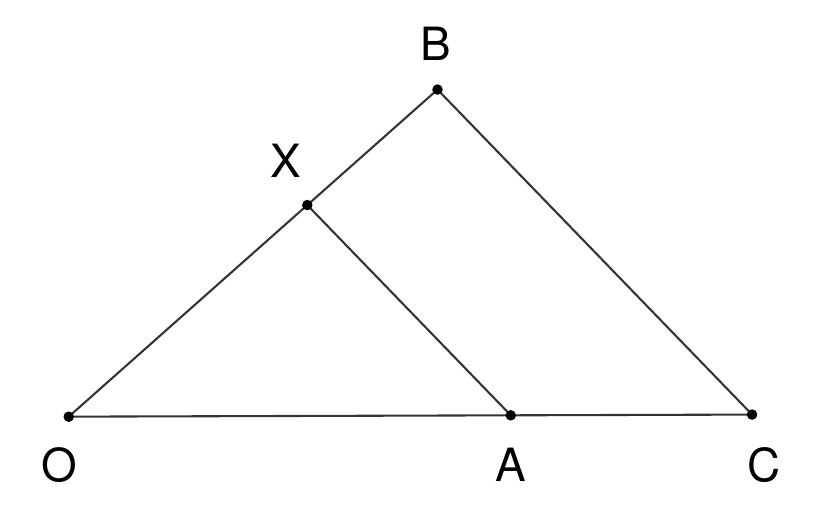}
\end{center}
\centerline {Figure 2}

\begin{rem}\label{r3}
Note that, given an interval $c$ and an angle $\al\in(0,\pi/2)$, one can construct a right triangle with the hypotenuse $c$ and an acute angle $\al$. Then, as the sides of this triangle are $a=c\sin\al$ and $b=c\cos\al$, we conclude that the intervals $a$ and $b$ are constructible as soon as $c$ and $\al$ are drawn.
\end{rem}

\section{Statement of the three-circle problem of Steiner}

\subsection{The Apollonius Three-Circle Problem}

Referring to Figure 3, we can state the classical Apollonius three-circle problem as follows:
\begin{prob} {\bf (Problem A)}\label{probA}
Given three circles $\mc C_1$, $\mc C_2$,  $\mc C_3$, construct a circle $\mc C$ that is tangent to each of them, that is, that intersects each at one point.
\end{prob}

\begin{center}
\includegraphics[width=0.5\textwidth]{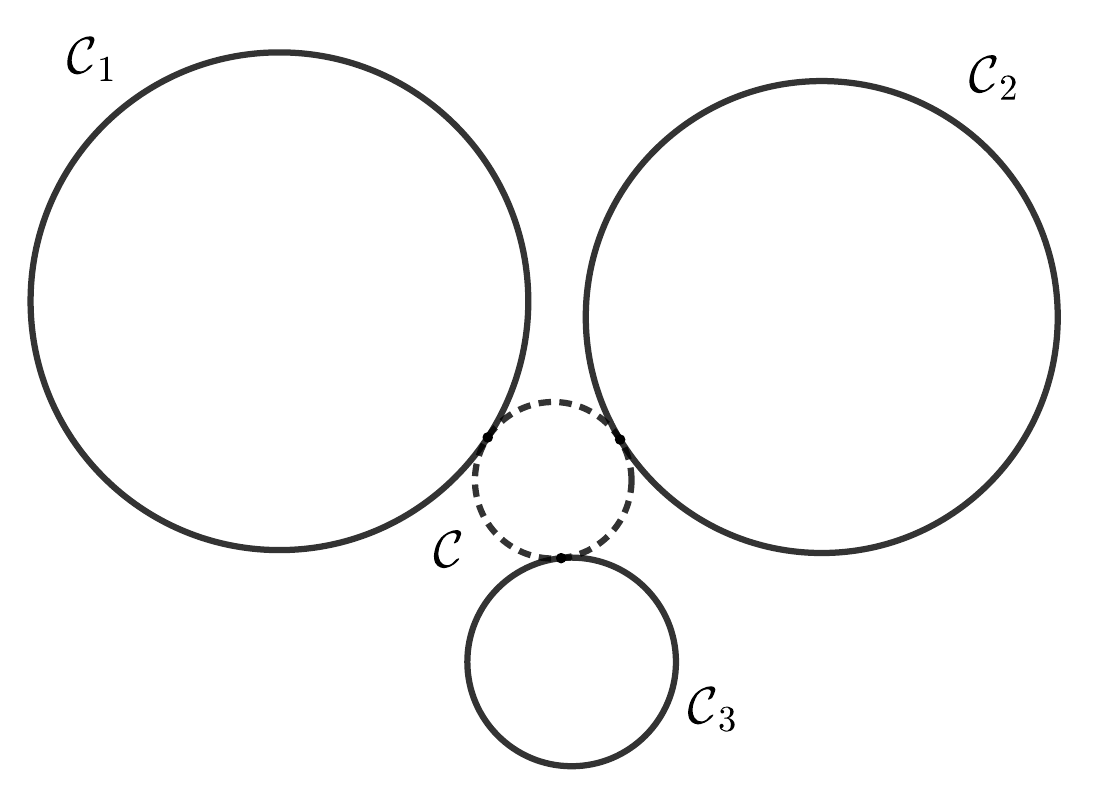}
\end{center}
\centerline{Figure 3}

\vskip 8pt
Note that, generically, there are eight different solutions to the Apollonius problem, as the circle $\mc C$ may touch one, or two, or all three given circles from outside.

\vskip 5pt
Figure 4 below is to state the following generalization of the Apollonius problem, that can be referred to as the three-circle problem of Steiner:

\begin{prob}\label{probS} {\bf (Problem S)}
Given three circles\, $\mc C_1$, $\mc C_2$, $\mc C_3$  with centers at $O_1$, $O_2$, $O_3$ and radii $r_1$, $r_2$, $r_3$, respectively, construct a circle $\mc C$ that intersects  $\mc C_1$ at an angle $0^0\leq\al\leq90^0$, $\mc C_2$ at an angle $0^0\leq\bt\leq90^0$, and $\mc C_3$ at an angle $0^0\leq\gm\leq90^0$.
\end{prob}

To avoid technicalities, we assume that the circles\, $\mc C_1, \mc C_2,\mc C_3$ are proper, their centers $O_1$, $O_2$, $O_3$ form a proper triangle, and that
\[
|O_1O_2|\ge r_1+r_2,\ \ |O_1O_3|\ge r_1+r_3, \text{ \ and\ } |O_2O_3|\ge r_2+r_3,
\]
as in Figure 4.

\vskip 8pt
\begin{center}
\includegraphics[width=0.65\textwidth]{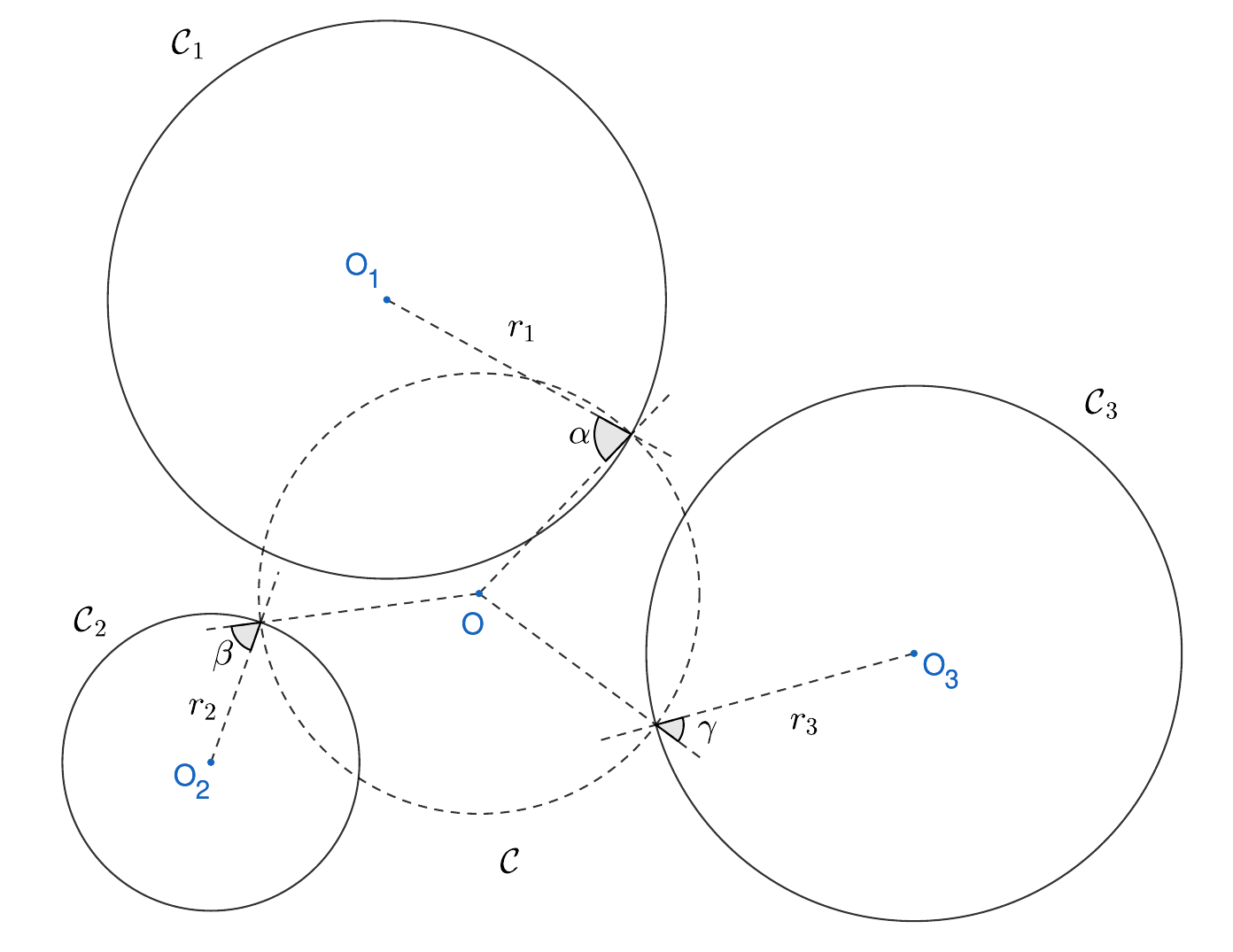}
\end{center}
\centerline {Figure 4}

\vskip 10pt
\section{Inversion and its properties}

\subsection{Definitions, basic properties and constructions}
Referring to Figure 5, consider a circle with center $O$ and radius $R$, which we will call the {\it inversion circle} and denote $\mc C_{\mc I}$.

\begin{center}
\includegraphics[width=0.5\textwidth]{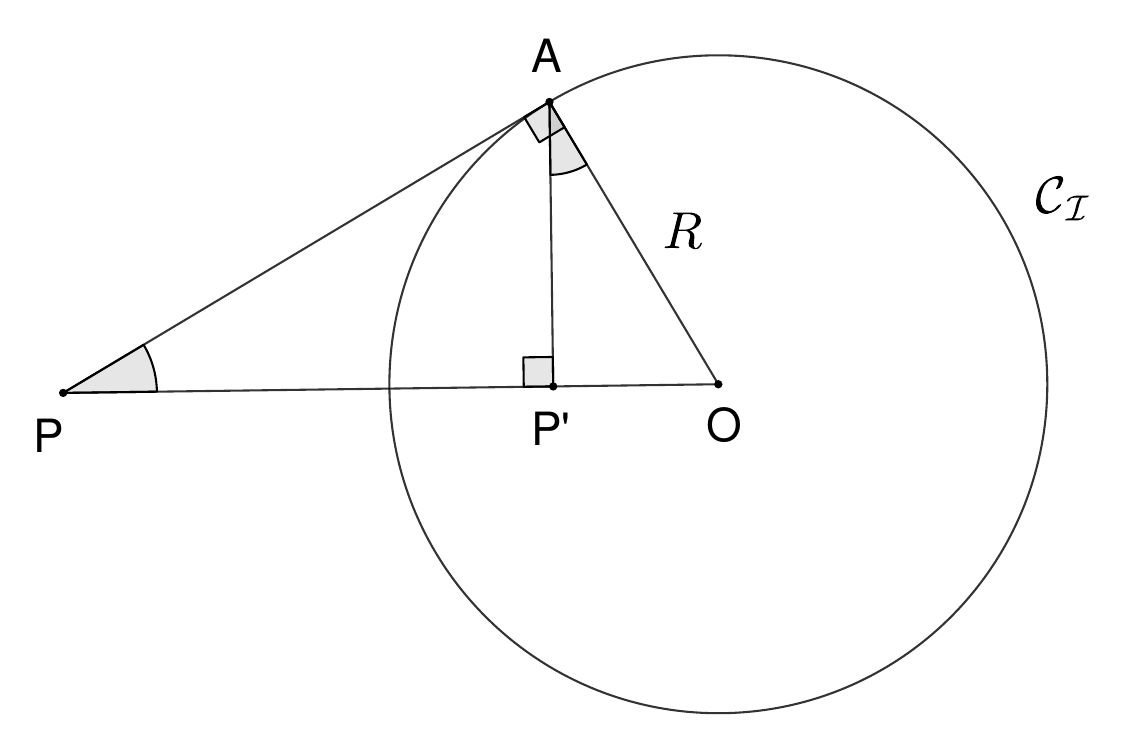}
\end{center}

\centerline{Figure 5}

\vskip10pt
\noindent
Points $P\ne O$ and $P^\pr\ne O$ are called the {\it inversion} images of each other with respect to $\mc C_{\mc I}$ if $P$ and $P^\pr$ lie on the same ray emanating from $O$ and
\[
|OP|\cdot|OP^\pr|=R^2.
\]

This transformation of the plane -- minus the point $O$ -- is called the {\it inversion} with respect to $\mc C_{\mc I}$. If we denote this inversion by $\mc I$, then  $\mc I(P)=P^\pr$, $\mc I(P^\pr)=P$, hence $\mc I^2(P)=P$ for every $P\ne O$, and $\mc I(P)=P$ whenever $P\in\mc C_{\mc I}$.

Figure 5 gives a way of constructing the point $P^\pr$ if $P$ is given and vise versa: If $P$ is given, let $A$ be the intersection of $\mc C_{\mc I}$ with the circle for which $[OP]$ is a diameter. Then $PA$ is tangent to  $\mc C_{\mc I}$, and if $AP^\pr\perp OP$, the triangles $AOP$ and $P^\pr OA$ are similar, so we have
\[
\frac{|OP|}R=\frac R{|OP^\pr|}\ \eqv \ |OP|\cdot|OP^\pr|=R^2.
\]
If $P^\pr$ is given, let $A\in\mc C_{\mc I}$ be such that $AP^\pr\perp OP^\pr$. Then, if $P\in OP^\pr$ is such that $AP\perp OA$, it follows as above that $|OP|\cdot|OP^\pr|=R^2$.

\vskip 5pt
Next, it can be proven -- see \cite[Theorem 5.45]{cg} -- that
\begin{enumerate}[(i)]
\item
If a circle $\mc C$ does not pass through the center $O$ of $\mc C_{\mc I}$, then $\mc I(\mc C)$ is a circle.
\item
If a circle $\mc C$ passes through $O$, then $\mc I(\mc C)$ is a straight line.
\end{enumerate}

\vskip 5pt
Let us now describe a construction that will need to perform repeatedly in the sequel without explicit reference:
\begin{prob}
Given a circle $\mc C$ that does not pass through the center of an inversion circle $\mc C_{\mc I}$, construct $\mc I(\mc C)$.
\end{prob}
\noindent
{\bf Solution:}
Pick three points $A, B, C$ on $\mc C$ and construct the points $\mc I (A)$, $\mc I (B)$, $\mc I(C)$ using the procedure given by Figure 5. Then construct the circle $\mc I(\mc C)$ as passing through the points $\mc I (A)$,
$\mc I (B)$, $\mc I(C)$ as in Problem \ref{prob1}.

\vskip5pt
Our solution of the three-circle problem of Steiner is based on the following facts:
\begin{enumerate}[(A)]
\item
An inversion preserves the angles between two circles -- see \cite[Theorem 5.51]{cg}. (Clearly, lines can be included in this statement, as circles with ``infinite radii.")
\item
Given two non-intersecting circles $\mc C_1$ and $\mc C_2$, there exists a circle $\mc C_{\mc I}$ such that $\mc C_1$ and $\mc C_2$ turn into concentric circles under the inversion with respect to $\mc C_{\mc I}$ -- see \cite[Theorem 5.71]{cg}.
\end{enumerate}

\subsection{Constructing a circle $\mc C_{\mc I}$ in (B)}
Referring to
Figure 6, we will describe the process of constructing possible centers for the circle $\mc C_{\mc I}$ from the statement (B) above as follows.

\begin{center}
\includegraphics[width=0.7\textwidth]{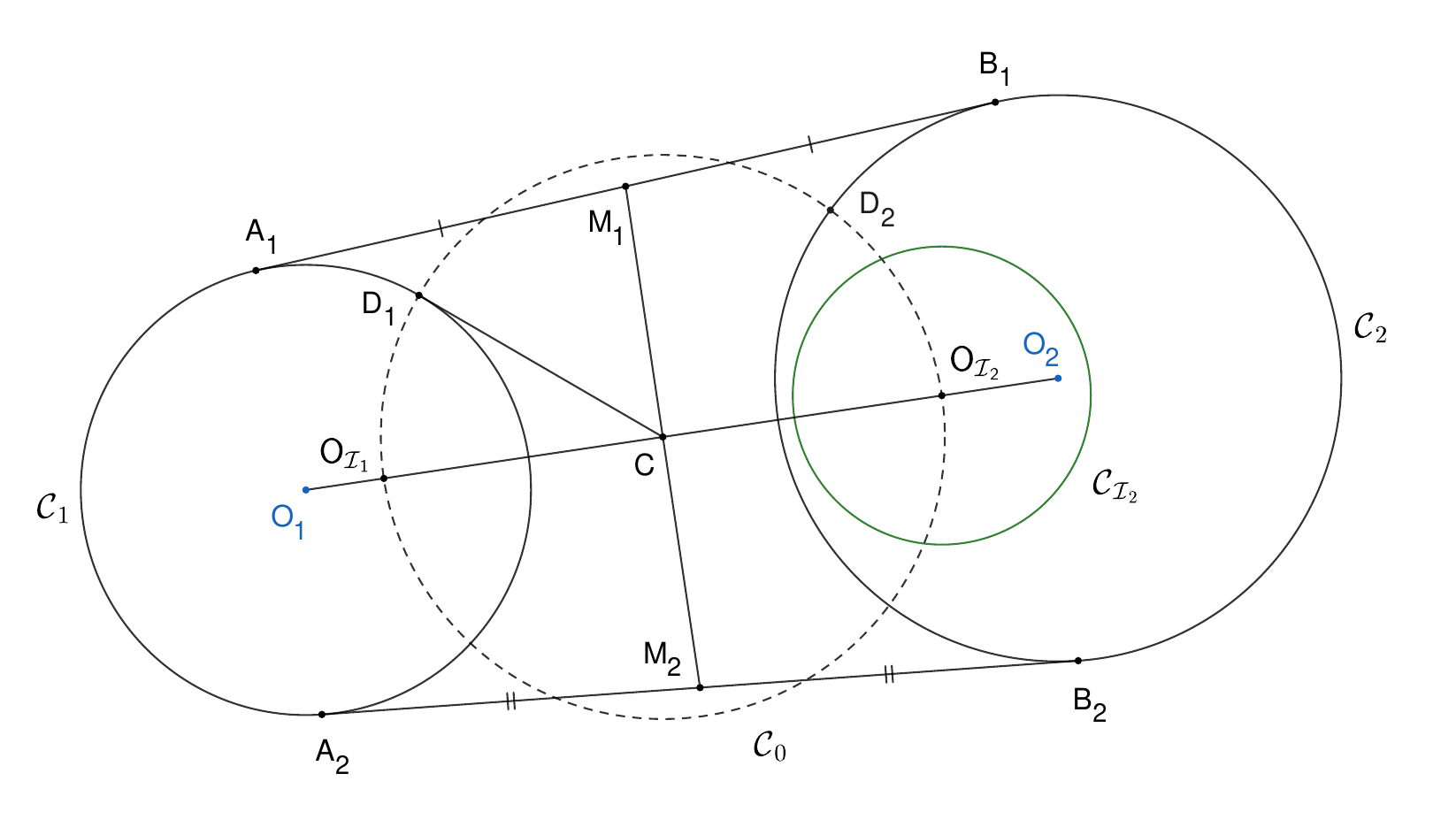}
\end{center}

\vspace{-5pt}
\centerline{Figure 6}

\vskip 8pt
Let $\mc C_1$ and $\mc C_2$ be non-intersecting circles with centers at $O_1$ and $O_2$ and $radii$ $r_1$ and $r_2$, respectively. Constructing a right triangle with the hypotenuse of length $|O_1O_2|$ and one of the sides of length $|r_1-r_2|$, we can draw the common outer tangents $[A_1B_1]$ and $[A_2B_2]$ to the circles $\mc C_1$ and $\mc C_2$, respectively.
Let $M_1$ and $M_2$ be the mid-points of segments $[A_1B_1]$ and $[A_2B_2]$, respectively.  It is clear that the line $M_1M_2$ is perpendicular to $O_1O_2$ and, following the argument in \cite[Theorem 2.21]{cg}, we conclude that  $M_1M_2$ is the {\it radical axis} of circles $\mc C_1$ and $\mc C_2$, that is, the tangent intervals to these circles eminating from any point on $M_1M_2$ are of the same length. Denote by $C$ the point of intersection of the lines $O_1O_2$ and $M_1M_2$.

Construct the point $D_1\in\mc C_1$ such that $[CD_1]$ is tangent to $\mc C_1$, as we did in Figure 5 with point $A$. Note that if $D_2$ is such a point on $\mc C_2$, then, since the point $C$ lies on the radical axis of $\mc C_1$ and $\mc C_2$, we have $|CD_2|=|CD_1|$.

Next, draw the circle of radius $|CD_1|$ with the center $C$ and denote it by $\mc C_0$. Since $CD_1$ and $CD_2$ are tangent to $\mc C_1$ and $\mc C_2$, respectively, $\mc C_0$ intersects these circles at $90$-degree angles.

Denote by $\mc O_{\mc I_1}$ and $\mc O_{\mc I_2}$ the points of intersection of $\mc C_0$ with $O_1O_2$, and let $\mc I$ stand for the inversion with respect to a circle of any non-zero radius centered at, say, $\mc O_{\mc I_2}$ (which can be clearly replaced by any circle centered at $\mc O_{\mc I_1}$). 
By (ii) in Subsection 4.1, $\mc I(\mc C_0)$ is a straight line, which, by (A) above, intersects the circles $\mc I(\mc C_1)$ and $\mc I(\mc C_2)$ at $90$-degree angles, hence passes through each of their centers.

On the other hand, since the center $\mc O_{\mc I_2}$ lies on the line $O_1O_2$, this line is invariant with respect to $\mc I$, and, as $O_1O_2$ is orthogonal to both $\mc C_1$ and $\mc C_2$, (A) above entails that $O_1O_2$ is orthogonal to $\mc I(\mc C_1)$ and $\mc I(\mc C_2)$, so, also passes through their centers.

Thus, we have two distinct lines $\mc I(\mc C_0)$ and $O_1O_2$, each passing through the centers of $\mc I(\mc C_1)$ and $\mc I(\mc C_2)$, implying that these centers coincide, that is, $\mc I(\mc C_1)$ and $\mc I(\mc C_2)$ are concentric.

\vspace{-2pt}
\section{Inversive geometry solution of the Problem of Steiner}

An elegant solution of the Apollonius three-circle problem -- in a particular case given by Figure 3 -- is given via an application of inversion as follows (see, for example \cite[pp.\,297-299]{ab}). First, the radii of the three given circles are enlarged by the same quantity so that two of the resulting circles become tangent at, say, a point $O$. Then an inversion with the center $O$ is  performed resulting in these two circles turn into parallel lines.
Note that property (B) above grants another striking approach to the Apollonius problem, one that uses a method of Gauss \cite{ga}.

\vskip5pt
\subsection{Investigation}
Assume first that there exists a circle $\mc C$ that satisfies the conditions in Problem S, as in Figure 4. Let $\mc I$ be the inversion from 4.2 with the inversion circle $\mc C_{\mc I}$ depicted in Figure 7 in green. Apply this inversion to the four other circles in Figure 7, the same as the four circles in Figure 4.

\vskip2pt

\begin{center}
\includegraphics[width=0.54\textwidth]{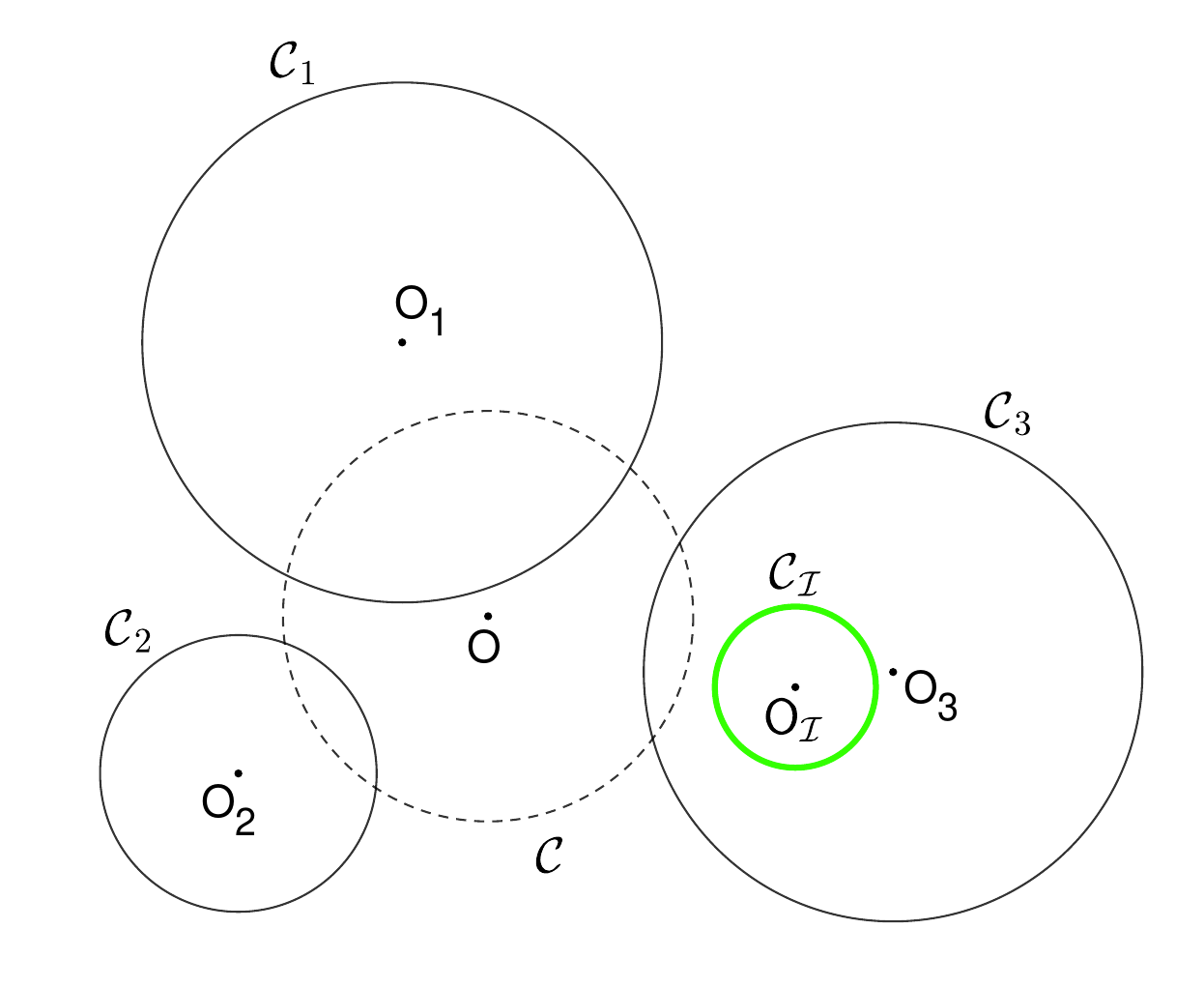}
\end{center}

\vspace{-5pt}
\centerline{Figure 7}

\vskip5pt
In view of the statement (A) above, the result looks as depicted in Figure 8, where\, $\mc C_1^\prime=\mc I(\mc C_1)$,\, $\mc C_2^\prime=\mc I(\mc C_2)$,\, $\mc C_3^\prime=\mc I(\mc C_3)$,\, $\mc C^\prime=\mc I(\mc C)$.
\vskip 10pt
\begin{center}
\includegraphics[width=0.75\textwidth]{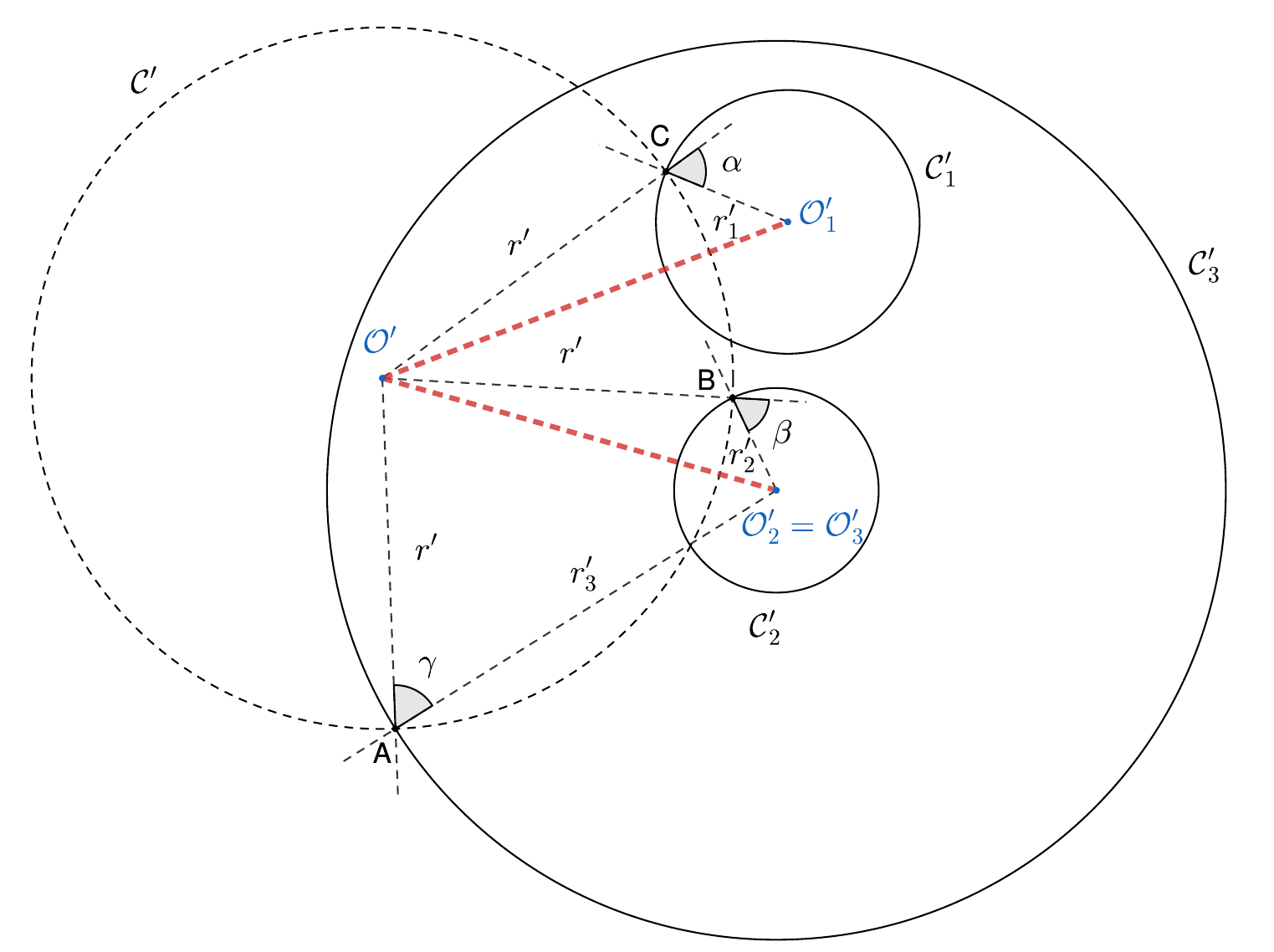}
\end{center}
\vskip3pt
\centerline{Figure 8}

\vskip 10pt
\noindent
Referring to Figure 8, we apply the Law of Cosines to triangles $\mc O^\pr A\mc O_2^\pr$ and $\mc O^\pr B\mc O_2^\pr$:
\[
(r^\pr)^2+(r_2^\pr)^2+2r^\pr r_2^\pr\cos\bt=(r^\pr)^2+(r_3^\pr)^2-2r^\pr r_3^\pr\cos\gm,
\]
implying that
\begin{equation}\label{e2}
2r^\pr(r_2^\pr\cos\bt+r_3^\pr\cos\gm)=(r_3^\pr)^2-(r_2^\pr)^2.
\end{equation}
We also calculate
\begin{equation}\label{e3}
|\mc O^\pr \mc O_1^\pr|=\sqrt{(r^\pr)^2+(r_1^\pr)^2+2r^\pr r_1^\pr\cos\al}
\end{equation}
and
\begin{equation}\label{e4}
|\mc O^\pr \mc O_2^\pr|=\sqrt{(r^\pr)^2+(r_2^\pr)^2+2r^\pr r_2^\pr\cos\bt}.
\end{equation}

\vskip 5pt
Now we are ready to solve Problem S:

\subsection{Solution of Problem S: constructing a circle $\mc C$}
Assume first that at least one of the angles $\al, \bt, \gm$ is less than $90^0$. Without loss of generality, let it be $\bt$ and/or $\gm$.  Then we can consider a circle $\mc C^\pr$ with radius
\begin{equation}\label{e8}
r^\pr=\frac{(r_3^\pr)^2-(r_2^\pr)^2}{2(r_2^\pr\cos\bt+r_3^\pr\cos\gm)}=\frac{r_3^\pr r_3^\pr}{2(r_2^\pr\cos\bt+r_3^\pr\cos\gm)}-\frac{r_2^\pr r_2^\pr}{2(r_2^\pr\cos\bt+r_3^\pr\cos\gm)}
\end{equation}
and center $\mc O^\pr$ such that $|\mc O^\pr \mc O_1^\pr|$ and $|\mc O^\pr \mc O_2^\pr|$ are given by (\ref{e3}) and (\ref{e4}), respectively, as depicted in Figure 9.

Note that, in view of Problem \ref{prob3}, and, as the interval $2(r_2^\pr\cos\bt+r_3^\pr\cos\gm)$, by Remark \ref{r3},  is constructible, the decomposition in (\ref{e8}) implies that we can construct the radius $r^\pr$. Next, construct a triangle with sides $r^\pr$ and $r_1^\pr$ and the angle $(180^0-\al)$ between them. By the Law of Cosines, the side opposite to this angle has length $|\mc O^\pr\mc O_1^\pr|$. Similarly, construct an interval of length $|\mc O^\pr\mc O_2^\pr|$. Then, with $|\mc O^\pr\mc O_1^\pr|$ and $|\mc O^\pr\mc O_2^\pr|$ constructed, we can consruct the point $\mc O^\pr$, hence, the circle $\mc C^\pr$.

\vskip 12pt
\begin{center}
\includegraphics[width=0.75\textwidth]{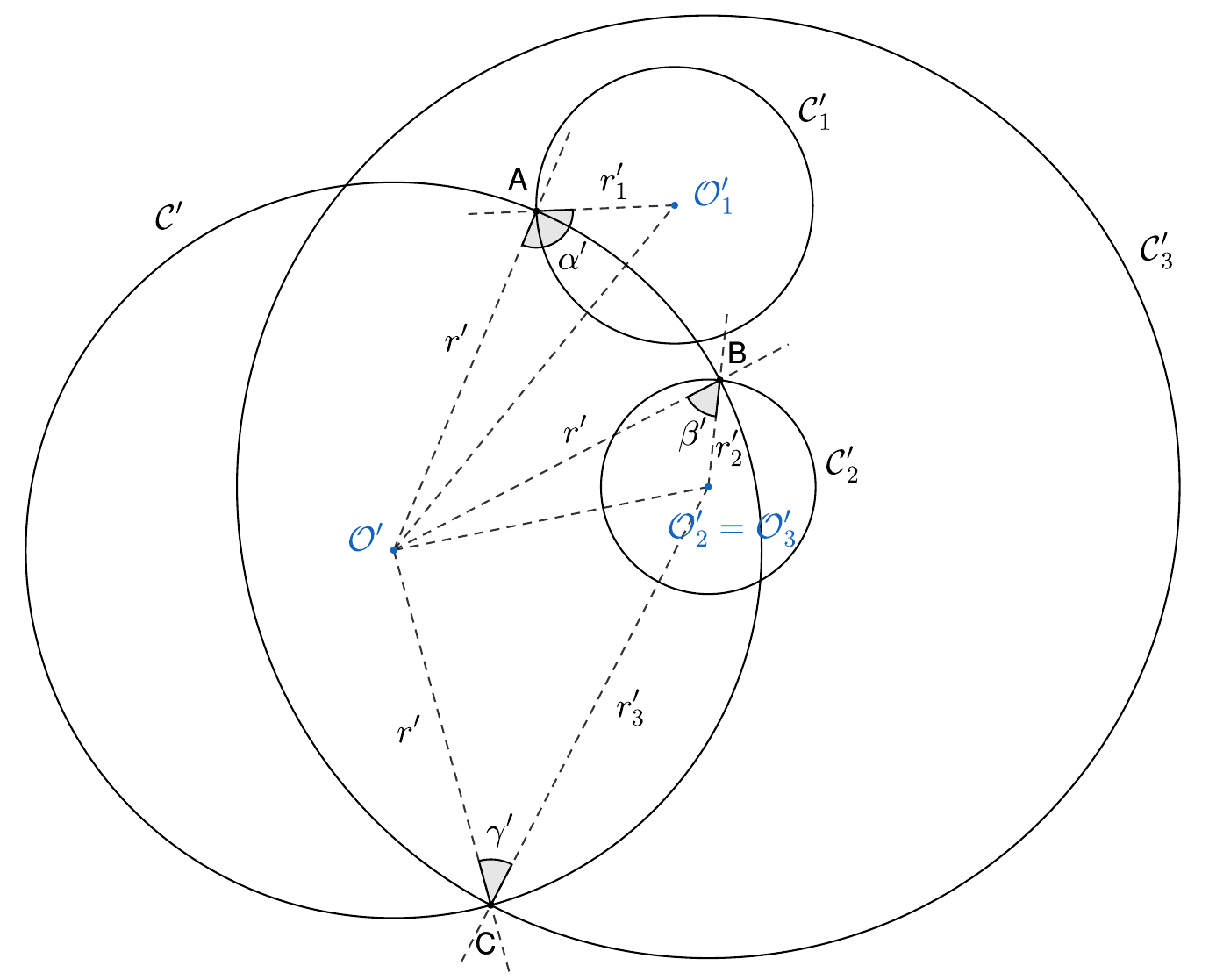}
\end{center}
\vskip 5pt
\centerline{Figure 9}

\vskip 10pt
Now, referring to Figure 9 and applying the Law of Cosines to triangle $\mc O^\pr A\mc O_1^\pr$, we obtain
\[
(r^\pr)^2+(r_1^\pr)^2-2r^\pr r_1^\pr\cos\al^\pr=(r^\pr)^2+(r_1^\pr)^2+2r^\pr r_1^\pr\cos\al,
\]
implying that
\begin{equation}\label{e5}
\cos\al^\pr=-\cos\al,\text{ \ hence \ }\al^\pr=180^0-\al.
\end{equation}
Similarly, considering triangle $\mc OB\mc O_2^\pr$, we infer
\begin{equation}\label{e6}
\cos\bt^\pr=-\cos\bt,\text{ \ hence \ }\bt^\pr=180^0-\bt.
\end{equation}
Finally, the Law of Cosines for triangle $\mc O^\pr C\mc O_2^\pr$ together with (\ref{e2}) yield
\[
2r^\pr(r_2^\pr\cos\bt+r_3^\pr\cos\gm^\pr)=(r_3^\pr)^2-(r_2^\pr)^2=2r^\pr(r_2^\pr\cos\bt+r_3^\pr\cos\gm),
\]
so
\begin{equation}\label{e7}
\cos\gm^\pr=\cos\gm,\text{ \ hence \ }\gm^\pr=\gm.
\end{equation}

Therefore, as (\ref{e5})\,-\,(\ref{e7}) indicate that the constructed circle $\mc C^\pr$ intersects $\mc C_1^\pr$ at the angle $\al$, $\mc C_2^\pr$ at the angle $\bt$, and $\mc C_3^\pr$ at the angle $\gm$, all that remains to do is to construct the circle $\mc C=\mc I(\mc C^\pr)$, which does satisfy conditions of Problem S.

\vskip 5pt
Now we are left with the case
\[
\al=\bt=\gm=90^0.
\]
\noindent
If we apply (construct) the inversion $\mc I$, then $\mc C^\pr=\mc I(\mc C)$ must intersect each of the concentric circles $\mc C_2^\pr=\mc I(\mc C_2)$ and $\mc C_3^\pr=\mc I(\mc C_3)$ at $90$-degree angles. This implies that $\mc I(\mc C)$ is a {\it line}, since there is no circle that is orthogonal to each of the two concentric circles -- in such a case, equation (\ref{e2}) would yield $r_2^\pr=r_3^\pr$, which is impossible.

\vskip5pt
\noindent
Note that, in view of the properties (i) and (ii) in Subsection 4.1, this implies that if a circle $\mc C$ intersects each of the two non-intersecting circles $\mc C_1$ and $\mc C_2$, as in Figure 6, at $90$-degree angles, then it passes through the center of $\mc C_{\mc I}$ (actually, it passes through both $O_1$ and $O_2$ in Figure 6). This is a striking fact by itself showcasing the remarkable usefulness of the inversion.

\vskip 5pt
\noindent
Observe that, being orthogonal to each of the circles\, $\mc C_1^\pr=\mc I(\mc C_1)$, $\mc C_2^\pr=\mc I(\mc C_2)$ (and $\mc C_3^\pr=\mc I(\mc C_3)$), the line $\mc I(\mc C)$ must pass through the centers $\mc O_1^\pr$ and $\mc O_2^\pr=\mc O_3^\pr$, so it can be drawn immediately. Constructing the inversion $\mc I$ of this line, we arrive at a circle $\mc C$ that intersects each of the circles\, $\mc C_1$, $\mc C_2$, $\mc C_3$ at a $90$-degree angle. This finalizes our solution of the three-circle problem of Steinter.

\begin{rem}
As with the classical Apollonius problem, there is more than one solution to Problem S. For example, one can see that changing $``-"$ for $``+"$ in front of $2r^\pr r_1^\pr\cos\al$ in (\ref{e3}) will not alter the values of $\al, \bt, \gm$. Therefore, such a change is admissible. The reader can find more admissible changes in definitions of $r^\pr$, $\mc O^\pr\mc O_1^\pr$, and  $\mc O^\pr\mc O_2^\pr$ resulting in different solutions of Problem S.
\end{rem}

\end{document}